\documentclass[11pt]{article}
\usepackage{amsmath}
\usepackage{mathrsfs}
\usepackage{amsthm}
\usepackage{amsfonts}
\usepackage{amssymb}
\usepackage{bm}
\usepackage{mathtools}
\usepackage{color}

\usepackage{enumerate}
\usepackage{enumitem}
\setlist[enumerate,1]{label=(\arabic*),font=\textup,
leftmargin=7mm,labelsep=1.5mm,topsep=0mm,itemsep=-0.8mm}
\setlist[enumerate,2]{label=(\alph*).,font=\textup,
leftmargin=7mm,labelsep=1.5mm,topsep=-0.8mm,itemsep=-0.8mm}

\newtheorem{theorem}{Theorem}[section]

\newtheorem{lemma}{Lemma}[section]

\title{\bf Domination in intersecting hypergraphs\thanks {Research was partially supported by the National Nature Science Foundation of China (No. 11571222)}}
\author {Yanxia Dong$^{1}$, \, Erfang  Shan$^{1,2}$\thanks{\em Corresponding authors. Email address: efshan@shu.edu.cn}, \, Shan Li, \, Liying Kang$^{1}$ \\
{\small $^{1}$Department of Mathematics, Shanghai University,
Shanghai 200444, P.R. China}\\
{\small$^{2}$School of Management, Shanghai University,
Shanghai 200444, P.R. China}}

\date{}
\begin{document}

\maketitle

\begin{abstract}

A matching in a hypergraph $H$ is a set of pairwise disjoint hyperedges. The matching number $\alpha'(H)$ of $H$ is the size of a maximum matching in $H$.  A subset $D$ of vertices of $H$ is  a dominating set of $H$ if for every $v\in V\setminus D$ there exists $u\in D$ such that $u$ and $v$ lie in an hyperedge of $H$. The cardinality of a minimum dominating set of $H$ is called the domination number of $H$, denoted by $\gamma(H)$.  It is known that for a intersecting hypergraph $H$ with rank $r$, $\gamma(H)\leq r-1$. In this paper we present structural properties on intersecting hypergraphs with rank $r$ satisfying the equality $\gamma(H)=r-1$. By applying the properties we show that all linear intersecting hypergraphs $H$ with rank $4$ satisfying  $\gamma(H)=r-1$ can be constructed by the well-known Fano plane.


\medskip

\noindent{\bf AMS (2000) subject classification:}  05C65, 05C69, 05C70
\end{abstract}

\section{Introduction}

Hypergraphs are a natural generalization of undirected graphs in which ``edges" may consist of more than 2 vertices.
More precisely, a (finite) {\em hypergraph} $H=(V(H), E(H))$ consists of a (finite) set $V(H)$ and a collection $E(H)$ of non-empty
subsets of $V(H)$. The elements of $V(H)$ are called {\em vertices} and the elements of $E(H)$ are called {\em hyperedges}, or simply {\em edges} of the
hypergraph.
An $r$-{\em edge} is an edge containing exactly $r$ vertices. The {\em rank} of  $H$, denoted by $r(H)$, is the maximum size of an edge in $H$. Specially, An $r$-{\em uniform} hypergraph $H$ is a hypergraph such that all edges are  $r$-edges.  A hypergraph is called {\em linear} if any two edges of the hypergraph intersect in at most one vertex. Obviously,   every (simple) graph is a linear 2-uniform  hypergarph.  Throughout this paper, we only consider hypergraphs of rank $r\geq2$ without multiple edges
and all edges have size at least $2$.

The {\em degree} of a vertex $v$ in $H$, denoted by $d_{H}(v)$, is the number of edges of $H$ containing the vertex $v$. A vertex
of degree zero is called an {\em isolated vertex}. A vertex of degree $k$ is called a {\em degree}-$k$ vertex.
The {\em minimum} and {\em maximum degree} among the vertices of $H$
are denoted by $\delta(H)$ and $\varDelta(H)$, respectively.
The {\em quasidegree} of  $v$ in $H$, denoted $qd_{H}(v)$,
is the maximum number of edges of $H$ whose pairwise intersection is only $v$.
Two vertices $u$ and $v$ in $H$ are {\em adjacent} if there is a edge $e$
of $H$ such that $\{u,v\}\subseteq e$. The {\em open neighborhood} of a vertex $v$ in $H$,
denoted $N_{H}(v)$, is the set of all
vertices different from $v$ that are adjacent to $v$.
If $H$ is clear from the context, we denote $d_{H}(v)$, $qd_{H}(v)$ and $N_{H}(v)$
 by $d(v), qd(v)$ and $N(v)$, respectively.
Two edges in $H$ are said to be {\em overlapping} if they intersect in at least two vertices.

A {\em partial hypergraph} $H'=(V(H'), E(H'))$ of $H=(V(H), E(H))$, denoted by $H'\subseteq H$, is a hypergraph such
that $V(H')\subseteq V(H)$ and $E(H')\subseteq E(H)$. In the class of graphs, partial hypergraphs are called {\em subgraphs}.
In particular, if $V(H')=V(H)$, $H'$ is called a {\em spanning partial hypergraph} of $H$.

For a hypergraph $H$ and  $X\subset V(H)$,  $H-X$
denotes the hypergraph obtained by removing the vertices $X$ from $H$ and removing
all edges that intersect $X$. For a subset $E'\subseteq E(H)$ of edges in $H$, we define $H-E'$ to be the hypergraph obtained from $H$ by deleting the edges in $E'$ and resulting isolated vertices, if any. If $E'=\{e\}$, then we write $H-E'$ simply as $H-e$. For $e\in E(H)$ and $v\in e$, if we remove the vertex $v$
from the edge $e$, we say that the resulting edge is obtained by $v$-{\em shrinking} the edge
$e$.

\subsection{Domination, matchings and transversals in hypergraphs}

A {\em dominating set} in a hypergraph $H$  is a subset $D$ of vertices of $H$  such that for every
vertex  $v\in V(H)\setminus D$ there exists an edge  $e\in E(H)$ for which  $v\in e$ and $e\cap D\neq \emptyset$.
Equivalently, every vertex
$v\in V(H)\setminus D$ is adjacent to a vertex in $D$. The minimum cardinality of a dominating set in $H$
is called its {\em domination number}, denoted by $\gamma(H)$.
A {\em matching} in $H$ is a set of disjoint edges. The maximum cardinality
of a matching in $H$ is called the {\em matching number}, denoted by $\alpha'(H)$.
A subset $T$ of vertices in $H$ is a {\em transversal} (also called {\em cover})
if $T$ has a nonempty intersection with each edge of
$H$. The {\em transversal number}, $\tau(H)$, is
the minimum size  of a transversal of $H$.
Transversals and matchings in hypergraphs are well studied
in the literature (see e.g. \cite{Alon,Chv,Dorf,fu,Henning1,Henning11,Henning2,Henning4}) and elsewhere. Domination in hypergraphs, was introduced  by Acharya \cite{Acharya1}
and studied further in  \cite{Acharya2,Aru,Bujt, Henning1,Jose,ShD}.

For a hypergraph $H$ of rank $r$, when $r=2$, $H$ is a graph,
Haynes et al. \cite{hhs1} observed that $\gamma(H)\leq \alpha'(H)\leq \tau(H)$.
When $r\ge 3$, by definitions, clearly $\gamma(H)\leq \tau(H)$ and  $\alpha'(H)\le \tau(H)$ still hold.
The extremal graphs, i.e., linear $2$-uniform hypergraphs achieving $\gamma(G)\leq \tau(G)$ were studied in \cite{Aru,lw,Ran1,Wu}
Recently, Arumugam et al. \cite{Aru} investigated the hypergraphs of rank $r\ge 3$ satisfying $\gamma(H)=\tau(H)$, and  proved that
their recognition problem is NP-hard  on the class of linear hypergraphs of rank 3.

In \cite{kld} we observed that the inequality $\gamma(H)\leq \alpha'(H)$  does not hold for a hypergraph $H$ of rank $r\ge 3$,
and the difference $\gamma(H)-\alpha'(H)$  can be arbitrarily large.
Further, we obtained  the following inequality.
\begin{theorem}(\cite{kld})\label{thm0}
If $H$ is a hypergraph  of rank $r\ge 2$ without isolated vertex, then $\gamma(H)\leq (r-1)\alpha'(H)$
 and this bound is sharp.
\end{theorem}
 In particular, if $r=2$ in Theorem \ref{thm0}, as Haynes et al. \cite{hhs1} observed,   $\gamma(H)\leq \alpha'(H)$.

For extremal hypergraphs of rank $r$ satisfying $\gamma(H)=(r-1)\alpha'(H)$,
Randerath et al. \cite{Ran2} gave a characterization of graphs (hypergraphs of rank $r=2$) with minimum degree two.  In 2010, Kano et al. \cite{Kano} provided a complete characterization of graphs with minimum degree one. For the case when rank $r=3$,  we give a complete characterization of hypergraphs $H$ in \cite{ShD}. For the case when $r\ge 4$, a constructive characterization of hypergraphs
with $\gamma(H)=(r-1)\alpha'(H)$ seems difficult to obtain. Thus we  restrict our attention to intersecting hypergraphs.

 A hypergraph is {\em intersecting} if any two edges have nonempty intersection. Clearly, $H$ is intersecting if and only if
$\alpha'(H)=1$. Intersecting hypergraphs are well studied in the literature (see, for example, \cite{ah,da, Ekr,fr, gui,henning5,hilton,man}).
For an intersecting hypergraph $H$ of rank $r$, we immediately have $\gamma(H)\leq r-1$.

In this paper  we first give some structural properties on the intersecting hypergraphs of rank $r$ achieving the equality $\gamma(H)=r-1$.
By applying the properties and Fano plane,
we provides a complete characterization of linear intersecting hypergraphs $H$ of rank $4$
satisfying $\gamma(H)=r-1$.

\section{The intersecting hypergraphs of rank $r$ with $\gamma(H)=r-1$}

In this section we give some  structural properties on intersecting hypergraphs of rank $r$ satisfying $\gamma(H)=r-1$.
The properties  play an important role in the characterization of intersecting hypergraphs of rank $4$ with $\gamma(H)=r-1$.

Let $\mathcal{H}_r$ be a family of intersecting hypergraphs of rank $r$
in which each hypergraph $H$ satisfies $\gamma (H)=r-1$.

\begin{lemma}\label{pro2.1}
 For every $H\in \mathcal{H}_r$, there exists an $r$-uniform spanning partial hypergraph $H^*$ of $H$ such that every edge in $H^*$ contains exactly one degree-$1$ vertex.
\end{lemma}

\noindent{\bf Proof.}
Let $H_0=H$. We define recursively the hypergraph $H_i$ by $H_{i-1}$. If there exists an edge $e_{i-1}\in E(H_{i-1})$ such that $d_{H_{i-1}}(v)\geq 2$ for each vertex $v$ in $e_{i-1}$, then set $H_{i}:=H_{i-1}-e_{i-1}$  for $i\ge 1$.
By repeating this process until every edge which remains contains at least one  degree-$1$ vertex, we obtain a spanning partial hypergraph of $H$.
Assume that the above process stops when $i=k$. Let $H^*=H_{k}$.  Then $H^*$ is a spanning partial hypergraph of $H$.
Clearly,   every edge in $H^*$ contains at least one degree-1 vertex and $H^*$ is still  intersecting.

 We claim that each edge in  $H^*$ contains exactly one degree-$1$ vertex.
 Suppose not. Then there exists an edge $e$ containing at least two degree-$1$ vertices.
 Let $D=\{v\in e ~|~d_{H^*}(v)\geq2\}$. Then $|D|\leq r-2$.
 Since $H^*$ is intersecting, $D$ is a transversal of $H^*$, so $\gamma(H^*)\leq\tau(H^*)\leq|D|\leq r-2$. Since $H^*$ is
a spanning partial hypergraph of $H$, we have $\gamma(H)\leq\gamma(H^*)\leq r-2$, contradicting  the assumption that $\gamma(H)=r-1$.
Further, we show that $H^*$ is $r$-uniform.  Suppose not. Let $e^{*}$ be an edge of $H^*$ such that $|e^{*}|\leq r-1$ and $u$  the unique degree-$1$ vertex of $e^{*}$. Since $H^*$ is intersecting, $e^{*}\backslash\{u\}$ is a dominating set of $H^*$.
Thus $\gamma(H)\leq\gamma(H^*)\leq r-2$, contradicting $\gamma(H)=r-1$ again.
\qed

For each $H\in \mathcal{H}_r$, let $H^*$ be the $r$-uniform spanning partial hypergraph of $H$  in Lemma \ref{pro2.1}.
Further, let $H'$ be the hypergraph obtained from $H^*$ by shrinking every edge to $(r-1)$-edge by removing the degree-1 vertex from each edge of $H^*$ and deleting multiple edges, if any. Obviously, $H'$ is an $(r-1)$-uniform intersecting hypergraph.

\begin{lemma}\label{pro2.2}
For every $H\in \mathcal{H}_r$, $\gamma(H)=\gamma(H^*)=\tau(H^*)=\tau(H')=r-1$.
\end{lemma}

\noindent{\bf Proof.}
Let $e\in E(H^*)$ and $u$ be the unique degree-1 vertex in $e$. Since $H^*$ is intersecting,
$e\setminus \{u\}$ is a transversal of $H^*$,  so $\gamma(H^*)\leq \tau(H^*)\leq r-1$.
On the other hand, note that $\gamma(H^*)\geq\gamma(H)=r-1$. Hence $\gamma(H^*)=\tau(H^*)=\gamma(H)=r-1$.
By the construction of $H'$, clearly any transversal of $H'$ is  a transversal of $H^*$
and $e\setminus\{u\}$ is also a transversal of $H'$. Hence $\tau(H')=\tau(H^*)$.
The equality chain follow. \qed

\begin{lemma}\label{pro2.3}
For $r\ge 3$ and every vertex $v$ in $H'$,  $2\leq qd_{H'}(v)\leq r-1$.
\end{lemma}

\noindent{\bf Proof.} Suppose, to the contrary, that there exists a vertex $v\in V(H')$ such that
$qd_{H'}(v)\le 1$ or $qd_{H'}(v)\ge r$.

Suppose that $qd_{H'}(v)\le 1$. Note that every vertex $H'$ has degree at least 2, so $qd_{H'}(v)\neq 1$.
Hence $qd_{H'}(v)=0$. Let $e$ is an edge containing $v$. Since $H'$ is intersecting, $e\cap f\neq \emptyset$ for any
$f\in E(H')\setminus \{e\}$. In particular, if $v\in e\cap f$, then $|e\cap f|\ge 2$ since $d_{H'}(v)\ge 2$ and $qd_{H'}(v)=0$.
Thus $e\setminus \{v\}$ would be a transversal of $H'$, contracting the fact in Lemma \ref{pro2.2}.

 Suppose that $qd_{H'}(v)\geq r$.
 Let $e_1,\ldots,e_{r}$ be the edges whose pairwise intersection is only $v$.
  By Lemma \ref{pro2.2} and $r\ge 3$,  we have $\tau(H')=r-1\ge 2$. This implies that
  there exists  an edge $g$ such that $v\not\in g$. Since $H'$ is  intersecting,
  $|g\cap e_i|\ge 1$ for each $i=1, 2,\ldots,r$. But then $|g|\ge r$,
   contracting the fact that $H'$ is $(r-1)$-uniform.
\qed
\begin{lemma}\label{pro2.4}
Let $H\in \mathcal{H}_r$ ($r\ge 3$). If $H$ is linear, then every edge of $H'$ has at most one degree-2 vertex and $\varDelta(H')=r-1$.
\end{lemma}
\noindent{\bf Proof.}
 If $H$ is linear, so is $H'$.
First, we show that every edge of $H'$ has at most one degree-2 vertex. Suppose not, and let $e\in E(H'),v_1,v_2\in e$ such that
$d_{H'}(v_1)=d_{H'}(v_2)=2$.  Then there exists  two distinct edges $e_1, e_2$ such that $e_i\cap e=\{v_i\}$ for $i=1,2$. Since $H'$ is a linear intersecting hypergraph, we have $|e_1\cap e_2|=1$, so there exists a vertex $x\in V(H')$ such that $e_1\cap e_2=\{x\}$. Then $T=(e\setminus \{v_1,v_2\})\cup \{x\}$ is a transversal of $H'$. Consequently, $\tau(H')\leq r-2$, a contradiction to $\tau(H')=r-1$.

  Next we show that $\varDelta(H')=r-1$.  Suppose not,  let $\varDelta(H')=r-a$ where $a\ge 2$.  As we have seen, $H'$ is a linear intersecting $(r-1)$-uniform hypergraph with $\tau(H')=r-1$.
Then $|V(H')|\le (r-1)(r-a)$. Let $v\in V(H')$ such that $d_{H'}(v)=r-a$. Then
$|N_{H'}(v)|=(r-2)(r-a)$. Note that $V(H')\setminus N_{H'}(v)$ is a transversal of $H'$.
But $|V(H')\setminus N_{H'}(v)|\le (r-1)(r-a)-(r-2)(r-a)=r-a\le r-2$. Thus $\tau(H')\le r-2$, contradicting
that $\tau(H')=r-1$.
\qed

\begin{lemma}\label{pro2.5}
Let $H\in \mathcal{H}_r$ ($r\ge 3$). If $H$ is linear, then  $3(r-2)\le |E(H')|\le (r-1)^2-(r-1)+1$, $n(H')=(r-1)^2-(r-1)+1$, and so
$\gamma(H')=1$.
\end{lemma}
\noindent{\bf Proof.}
 Since $H'$ is a linear intersecting $(r-1)$-uniform hypergraph, $|E(H')|=\sum_{v\in e}d_{H'}(v)-(r-2)$ for any edge
 $e\in E(H')$. By Lemma \ref{pro2.4}, we immediately have $3(r-2)\le|E(H')|\le (r-1)^2-(r-1)+1$.

 We now show that $n(H')=(r-1)^2-(r-1)+1$. Let $v\in V(H')$ such that $d_{H'}(v)=\Delta(H')=r-1$. Then
  $|n(H')|\geq |N_{H'}(v)\cup\{v\}|=(r-1)(r-2)+1=(r-1)^2-(r-1)+1$. Suppose that $|n(H')|\geq (r-1)^2-(r-1)+2$. Then there exists $u\in V(H')$ such that
  $u\not\in N_{H'}(v)\cup\{v\}$.
By Lemma \ref{pro2.3}, $d_{H'}(u)\ge 2$, so there exist two edges $e_1$ and $e_2$ such that $u\in e_1 \cap e_2$.
Clearly, $v\not\in e_1$ and $v\not\in e_2$. Since $H'$ is linear intersecting, $e_1$ intersects each one of the edges that contains $v$, implying that $|e_1|\geq 5$. This contradicts that $H'$ is an $(r-1)$-uniform hypergraph.  Therefore, $n(H')=(r-1)^2-(r-1)+1$, that is, $n(H')=|N_{H'}(v)\cup\{v\}|$.
This implies that $\{v\}$ is a dominating set of $H'$, so $\gamma(H')=1$.
\qed

\section{Linear intersecting hypergraphs $H$ of rank $4$ with $\gamma(H)=3$}

In the section we give a complete characterization of linear intersecting hypergraphs $H$ of rank $4$ with $\gamma(H)=r-1$.
For this purpose, let $F$  be the Fano Plane and let $F^{-}$ be the hypergraph obtained from $F$ by deleting any edge of $F$.
The two hypergraphs $F$ and $F^-$ are shown in Fig. 1.
\begin{figure}[htbp]
  \centering
  \includegraphics[height=4.2cm]{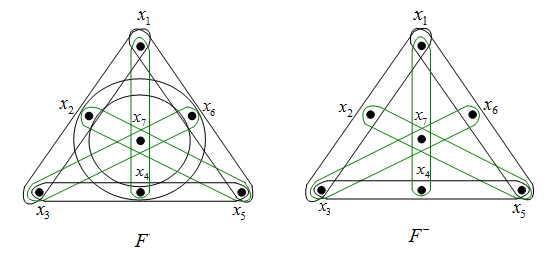}
  \caption{The Fano Plane $F$ and the hypergraph $F^{-}$}\label{Fig.1}
\end{figure}

\begin{lemma}\label{thm4.1}
let $H\in \mathcal{H}_4$ and  $H'$ be the hypergraph as defined in the above section. If $H$ is linear,  then
$H'= F$ or $H'= F^{-}$.
\end{lemma}

\noindent{\bf Proof.}
By Lemma \ref{pro2.5}, we have $6\leq |E(H')| \leq 7$ and $n(H')=7$ for $r=4$.
Note that $H'$ is a linear intersecting 3-uniform hypergraph. If $|E(H')|=7$, then $H'$ must be the Fano plane $F$.
If $|E(H')|=6$, then $H'$ is the hypergraph $F^-$ (see Fig. 1).
\qed

To complete our characterization, we let $F_1$ ($F_1^-$) be the hypergraph obtained from $F$ ($F^-$) by adding a new vertex
to each edge of $F$  ($F^-$), respectively. Let $F_2$ be the hypergraph obtained from $F_1$ by shrinking one edge
to $3$-edge by removing the degree-1 vertex in the edge. Let $F_3$  be the hypergraph obtained from $F_1^-$ by adding a new edge
$f=\{x_{i_1}, x_{i_2}, x_{i_3}, x_{i_4}\}$ where $x_{i_1}, x_{i_2}, x_{i_3}$ and  $x_{i_4}$ lie
in distinct edges of $F_1^-$ and $d_{F_1^-}(x_{i_1})=d_{F_1^-}(x_{i_2})=2$, $d_{F_1^-}(x_{i_3})=d_{F_1^-}(x_{i_4})=1$.
We define $\mathcal{L}=\{F_1, F_1^-, F_2, F_3\}$ (see Fig. 2). Clearly, every hypergraph in $\mathcal{L}$ is a linear intersecting hypergraph of rank 4.

\begin{figure}[htbp]\label{fig2}
  \centering
  \includegraphics[height=4.2cm]{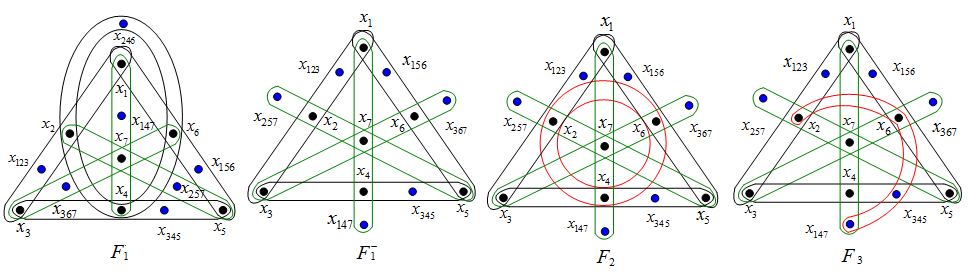}
  \caption{The hypergraphs $F_1, F^-_1, F_2$ and $F_3$}\label{Fig.2}
\end{figure}

\begin{theorem}\label{thm4.2}
For a linear intersecting hypergraph $H$ of rank $4$,  $\gamma(H)=3$ if and only if
$H\in \mathcal{L}$.
\end{theorem}
\noindent{\bf Proof.}
First, suppose that $H\in \mathcal{L}$, and let $e$ be an arbitrary edge of $H$ containing four vertices and $v$
the degree-1 vertex.
Then it is easy to check that $D=e\setminus \{e\}$ is a minimum dominating set of $H$. Thus $\gamma(H)=3$.

Conversely, suppose that $\gamma(H)=3$, we show that $H\in \mathcal{L}$.
Let $H^*$ and $H'$ be the hypergraphs corresponding to $H$ as defined in above section.
By Lemma \ref{thm4.1}, $H'=F$
or $H'=F^{-}$, so $H^*= F_{1}$ or $H^*= F^{-}_{1}$.

 {\it Case} 1. $H^*=F_{1}$. In this case, we claim that $H=H^*=F_1$. It suffices to show that
$E(H)\setminus E(H^*)=\emptyset$.
Suppose not. Let $e\in E(H)\setminus E(F_1)$. Then $|e\cap f|=1$ for any $f\in E(F_1)$.
By the construction of $F_1$,  $V(F_1)=V(F)\cup I$ where $I$ consists of seven degree-1 vertices in $F_1$. Note that any two vertices of $F$ lie in exactly one common edge of $F_1$, so $|e\cap V(F)|\le 1$. This implies that $|e|\ge 5$, since $H$ is linear and intersecting.
This contradicts that $r(H)=4$.

 {\it Case} 2. $H^*= F^{-}_{1}$. In this case, we show that $H\in \{F^-_1, F_2, F_3\}$.
 It suffices to show that $H=F_3$ if $H\neq F_1^-$ and $H\neq F_2$.
Let $V(F^-_1)=V_1\cup V_2\cup V_3$ where $V_i$ is the set of degree-$i$ vertices in $F^-_1$. Then $|V_1|=6, |V_2|=3$ and $|V_3|=4$.
 Suppose now that $H\neq F_1^-$ and $H\neq F_2$.
 Then $E(H)\setminus E(F^-_1)\neq\emptyset$. Let $e\in E(H)\setminus E(F^-_1)$.
 Suppose that $e\cap V_3\neq \emptyset$. Since $H$ is linear and intersecting, $|e\cap f|=1$ for any $f\in E(F^-_1)$.
 Note that any two vertices of $V_3$ lie in exactly one common edge of $F^-_1$, so $|e\cap V_3|=1$. Let $e\cap V_3=\{x_i\}$.
 This implies that $e\supseteq V_1\cap V(H-x_i)$. But then $x_i$ is a dominating set of $H$, contradicting that $\gamma(H)=3$.
 Hence $e\cap V_3=\emptyset$, and thus $e\subseteq V_1\cup V_2$. Suppose that $|e\cap V_2|\le 1$. Then $|e\cap V_1|\ge 4$.
 Hence $|e|\ge 5$, a contradiction. So $|e\cap V_2|=2$ since $H\neq F_2$. It immediately follows that $E(H)\setminus E(F^-_1)=\{e\}$.
 Therefore, $H=F_3$.
\qed

\section{Conclusions}

In this paper
we present the propositions of the intersecting hypergraphs  that achieve the equality $\gamma(H)=r-1$.
Especially, we provide a complete characterization of the linear intersecting hypergraphs with rank $r=4$ satisfying $\gamma(H)=3$.
One  is interested in characterizing the extremal intersecting hypergraphs with rank $r=4$ satisfying $\gamma(H)=3$.

\end{document}